\documentclass{amsart}
\usepackage{amsmath, amssymb, amsthm}
\usepackage{caption}
\usepackage{graphicx}
\usepackage{tabularx}
\usepackage{url}

\theoremstyle{plain}

\theoremstyle{definition}
\newtheorem*{definition}{Definition}

\newcommand{\vol}[1]{\operatorname{vol}(#1)}
\newcommand{\Cov}[1]{\operatorname{Cov}\left({#1}\right)}
\DeclareMathOperator*{\argmin}{arg\,min}

\begin{document}

\title[Deep Learning in Knot Theory]{Narrowing the Gap between Combinatorial and Hyperbolic Knot Invariants via Deep Learning}

\author{Daniel Grünbaum}

\address{Fakult\"at f\"ur Mathematik, Universit\"at Regensburg, 93040 Regensburg}
\email{daniel.gruenbaum@ur.de}

\begin{abstract}
We present a statistical approach for the discovery of relationships between mathematical entities that is based on linear regression and deep learning with fully connected artificial neural networks. The strategy is applied to computational knot data and empirical connections between combinatorial and hyperbolic knot invariants are revealed.
\end{abstract}

\maketitle

\section{Introduction}
\label{cha:intro}
A central goal in mathematical research is the discovery of theoretical relationships between different entities. After a relationship is noticed, mathematicians try to formalize all involved concepts by definitions, the exact nature of the relationship by theorems and the precise reasoning for the existence of the relationship by proofs. The first step of finding possible connections that warrant further investigation, however, is often left to a mixture of the researcher's intuition and domain knowledge.

In disciplines other than mathematics, for example in biology or in the social sciences, a routinely
performed exploratory step is gathering data for different entities and plotting them against each
other. In the case of an existing relationship, a visual inspection may reveal a non-random pattern
to the researcher and the strength of the association can subsequently be made more precise through
statistical evaluations. For a mathematician following this approach, several obstacles are bound to
occur: Firstly, data collection is traditionally not a prominent part of pure mathematics. Many
entities are hard to calculate by hand and computational aid in the form of specialized software is
not always available. Secondly, mathematical objects can be of infinite or at least prohibitively
high dimension, rendering a visualization impossible. Thirdly, a quantitative description of the
strength of a relationship through statistical methods does not constitute a formal theorem, let
alone a proof that shows the exact nature of the connection in rigorous mathematical terms. The goal
of this article is to argue that a data-driven exploratory step can still be a valuable tool for
knot theorists. With respect to the first of the above obstacles, scarcity of data, the field of
knot theory fares better than most other areas of mathematics due to an active community of
computational knot theorists that provide databases containing many knot invariants or suitable
tools for calculating them \cite{SnapPy, sage, khoho, knotinfo}. Given the amount of data that can
be generated, data-hungry algorithms of high-dimensional statistical learning, with deep learning as
its most prominent subfield \cite{goodfellow}, become viable options for replacing the visual
inspection. The last objection still holds, but the presented investigation of different knot
invariants will show that methods of data analysis can at least provide concrete research questions,
which can then be addressed by more traditional mathematical methods.

This article builds on the work of other researchers. A strong correlation between the determinant
and the hyperbolic volume of a knot has been established by Dunfield through visual inspection in an
unpublished note \cite{dunfield_determinant}.
The subject was revisited by Friedl and Jackson
\cite{friedl_hyperbolic}, who provided a quantitative analysis of the relationship between the
hyperbolic volume and several invariants derived from the Alexander polynomial, employing linear
regression as a statistical method.
Deep learning has been introduced to knot theory only recently:
Jejjala, Kar and Parrikar \cite{dl_jejjala} reported that a standard version of an artificial neural
network could predict the hyperbolic volume of a knot with high precision from the coefficients of
its Jones polynomial. They suggested incorporating invariants from Khovanov homology for a possible
improvement of the predictive power.
Following up on this discovery, Craven, Jejjala and Kar
\cite{craven} used layer-wise relevance propagation to further analyze and simplify the complex
neural network. They were able to distill a performant approximate formula for the hyperbolic volume that
inputs only a single evaluation of the Jones polynomial at a root of unity.
Hughes \cite{hughes} successfully used deep learning techniques to
predict from a braid word representation whether a given knot is quasipositive.
Joining forces, Craven, Hughes, Jejjala and Kar \cite{craven2} used a similar methodology to unveil empirical
connections between the Jones polynomial, Khovanov homology, smooth slice genus and
Rasmussen's $s$-invariant.
In an independent stream of research, Davies et al. \cite{davies} combined deep learning with
attribution techniques to detect a simple relationship between the knot signature, hyperbolic
volume, meridian and longitude translation. The empirical result lead to the formulation of a
conjecture that could be proven rigorously by adding the injectivity radius as a last ingredient \cite{davies_signature}.

In this article, we will use techniques from deep learning as an exploratory tool to gather evidence
for the existence of theoretical relationships between combinatorial and hyperbolic knot invariants.
The experimental evidence will by no means obliterate the necessity for a formal mathematical
treatment of the underlying abstract causes. On the contrary, it indicates only a starting point for
mathematical research of a more traditional nature.

The article is structured as follows:
Section~2 presents linear regression and deep learning as empirical methods of scientific discovery.
Section~3 briefly introduces combinatorial and hyperbolic knot invariants. The conducted experiments
and their results are documented and discussed in Section~4. Section~5 uses the discovered empirical
relationships to formulate a list of questions for subsequent theoretical research.

\section{Statistical Methods}
This section explains basic concepts of linear regression and deep learning, before arguing how prediction errors of artificial neural networks can be interpreted as a more evolved version of visual inspection methods and correlation coefficients.

\subsection{Linear Regression}
If we want to investigate the existence of a relationship between two $1$-dimensional entities, the
easiest way is to plot them against each other and inspect the resulting data cloud visually. An
example would be the aforementioned plotting of the hyperbolic volume against the knot determinant
for different classes of knots. If a relationship exists, the data shows a systematic pattern that
can be recognized by the human eye. In several fields of science, such as biology or the social
sciences, it is customary to assess the strength of a conjectured linear relationship quantitatively
via the Pearson correlation coefficient.

\begin{definition}
	Let $x= (x_i)_{i=1,...,n}$ be a finite family of real numbers.
	\begin{itemize}
		\item The \emph{sample mean} of $x$ is given by
		$\bar{x} = \frac{1}{n}\sum_{i=1}^{n} x_i$.
		\item The \emph{sample variance} of $x$ is given by
		$\sigma^2_x = \frac{1}{n-1}\sum_{i=1}^{n} (x_i-\bar{x})^2$.
		\item The \emph{sample standard deviation} of $x$ is given by
		$\sigma_x = \sqrt{\sigma^2_x}$.
		\item Let $y = (y_i)_{i=1,...,n}$ be another finite family of real numbers. The \emph{sample covariance} of $x$ and $y$ is given by
		$\Cov{x,y} = \frac{1}{n-1}\sum_{i=1}^{n} (x_i-\bar{x})(y_i-\bar{y})$.
		\item For $\sigma_x, \sigma_y \neq 0$, the \emph{Pearson correlation coefficient} of $x$ and $y$ is given by
		\begin{equation}
			r_{x,y} = \frac{\Cov{x,y}}{\sigma_x \sigma_y}.
		\end{equation}
	\end{itemize}
\end{definition}

Note that we have $r_{x,y} \in [-1,1]$ and $r_{x,y} \in \{-1,1\}$ if and only if $y$ is an affine
function of $x$. Visually, $\vert r_{x,y} \vert$ describes how well the data points $(x_i,y_i)$ fit
on a straight line. A related way of quantifying the strength of the association between $x$ and $y$ is
\emph{linear regression}, i.e. fitting the optimal linear approximation
\begin{equation}
	y \approx f(x) = a_* x + b_*
\end{equation}
to the data according to

\begin{equation}
	(a_*,b_*) = \argmin_{(a,b)\in \mathbb{R}^2} \frac{1}{n}\sum_{i=1}^{n} \Vert y_i - (a x_i+b) \Vert^2
\end{equation}
and using the expression in the minimum, which is called the \emph{mean squared error} $\Delta_{\text{MSE}}$, as a figure of merit. An additional relative error measurement is given by the \emph{mean absolute percentage error}

\begin{equation}
	\Delta_{\text{MAPE}}(x,y,a_*,b_*) = \frac{1}{n}\sum_{i=1}^{n} \frac{|y_i - (a_* x_i + b_*)|}{|y_i|}
\end{equation}
if $y$ has only non-zero entries.
Two main adjustments to this procedure are possible:
\begin{itemize}
	\item If the visual data inspection shows a non-linear pattern, we can transform $x$ by applying
	a suitable function $t \colon \mathbb{R} \to \mathbb{R}$ in such a way that the relationship
	between $t(x)$ and $y$ is linear. After the transformation, the same measures of dependence as
	above can still be used for an evaluation. A common example for this technique is the use of
	logarithmic scaling to quantify the strength of exponential relationships.
	\item If $x \in \mathbb{R}^m$ is a multidimensional entity and we suspect a linear relationship, we can simply fit a multilinear model $y \approx \beta \cdot x + b$ for $\beta \in \mathbb{R}^m$ that minimizes the mean squared error as above. Note that the ideal coefficients can be found analytically without any numerical optimization techniques.\footnote{Extensions to the case where $y$ is also multidimensional can be formulated in a similar way, but it is not necessary for our purposes.}
\end{itemize}

\subsection{Deep Learning}
The above adjustments are not sufficient if we want to investigate relationships that are non-linear
and involve high-dimensional data at the same time. Since we are looking for the existence of any,
possibly non-linear, relationship between $x$ and $y$, it is not sufficient to calculate the Pearson
correlation between $x$ and $y$ and give up if the absolute value is low. Contrary to the
$1$-dimensional case, we cannot simply plot the data and identify an ideal transformation $T \colon
\mathbb{R}^m \to \mathbb{R}^k$ for some $k \in \mathbb{N}$ that ensures a linear relationship
between $T(x)$ and $y$. Deep learning circumvents this problem by choosing a wide family of possible
parameterized transformations and optimizing not only the parameters of the multilinear model but
also those of the transformations. If the family of allowed transformations provides enough
expressive power and a suitable optimization strategy is chosen, we can be confident that any
possible functional relationship between $x$ and $y$ will be detected \cite{function_approximators}.
A common and well-studied form of deep learning is via fully connected artificial neural networks,
which we will define after a short notational digression.

\begin{definition}
	Let $k,l \in \mathbb{N}$. For a map $\phi \colon \mathbb{R}^k \to \mathbb{R}^l$ and a function $\psi\colon \mathbb{R} \to \mathbb{R}$, we define the \emph{componentwise composition of $\phi$ and $\psi$} to be the map
	\begin{equation}
		\begin{split}
		\psi \odot \phi \colon \mathbb{R}^k &\to \mathbb{R}^l \\
		x &\mapsto \left(\psi(\phi(x)_1),...,\psi(\phi(x)_l)\right).
		\end{split}
	\end{equation}
\end{definition}

The componentwise composition allows us to formulate the definition of an artificial neural network in a compact way:

\begin{definition}
	Let $m,p \in \mathbb{N}$ and let $(k_0,...,k_{p+1})$ be a family of natural numbers with $k_0=m$ and $k_{p+1} = 1$.
	\begin{itemize}
		\item A \emph{fully connected artificial neural network (ANN) with $p$ hidden layers} is a function
		$f \colon \mathbb{R}^m \to \mathbb{R}$
		that can be written as a composition
		\begin{equation}
			f = \phi_p \circ (\psi_{p-1} \odot \phi_{p-1}) \circ ... \circ (\psi_1 \odot \phi_1) \circ (\psi_0 \odot \phi_0)
		\end{equation}
		for a family of affine maps
		\begin{equation}
			\left(\phi_t \colon \mathbb{R}^{k_t} \to \mathbb{R}^{k_{t+1}}\right)_{t \in \{0,...,p\}}
		\end{equation}
		and a family of functions
		\begin{equation}
			\left(\psi_t \colon \mathbb{R} \to \mathbb{R}\right)_{t \in \{0,...,p-1\}}.
		\end{equation}
		For $t=1,...,p$, we call $k_t$ the \emph{number of hidden neurons in the $t$-th layer}.
		\item For each $t \in \{0,...,p\}$, there is a pair $(W_t, b_t)$
		of a $(k_{t+1} \times k_t)$-dimensional matrix $W_t$ with real entries and a vector $b_t \in \mathbb{R}^{k_{t+1}}$ such that we have
		\begin{equation}
			\phi_t(x) = W_tx +b_t
		\end{equation}
		for all $x \in \mathbb{R}^{k_t}$. We call $W_t$ the \emph{$t$-th weight matrix} and $b_t$ the \emph{$t$-th bias vector}. The entries of the weight matrices are called \emph{weights}.
		\item For each $t \in \{0,...,p-1\}$, we call $\psi_t$ the \emph{$t$-th activation function}.
	\end{itemize}
\end{definition}

As in the case of a multilinear model, we can view the entries of the weight matrices and bias vectors as free model parameters of the artificial neural network and obtain a parameterized model of the form
\begin{equation}
	\mathbb{R}^m \times \mathbb{R}^{(k_0 k_1) +(k_1 k_2) + ... + (k_p k_{p+1})} \times \mathbb{R}^{k_1 + ... + k_{p+1}} \to \mathbb{R}
\end{equation}
whose mean squared error on the data can be minimized with respect to the model parameters. Finding
the optimal pair $(W_p,b_p)$ for the last transformation $\phi_p$ is precisely a multilinear
regression. Since we optimize not only the weights for the regression function, but also the weights
in the previous matrices, the goal of optimizing the parameters of an ANN can be formulated as
follows: We want to find a transformation $T \colon \mathbb{R}^m \to \mathbb{R}^{k_p}$ that is
optimal for a prediction of the target variable~$y$ by multilinear regression and that can be
obtained by the application of affine transformations and activation functions.

The activation
functions enable non-linear transformations. Historically, \emph{sigmoidal} activation functions,
such as the \emph{hyperbolic tangent} $\psi(z) = \tanh(z)$ and the \emph{logistic function} $\psi(z)
= \frac{1}{1+e^{-z}}$ have been popular \cite{activation_glorot}. Named after the characteristic S-shape of their graphs,
these functions are monotonic and bounded. Their use is justified from a biological viewpoint where
neurons in the artificial neural network correspond to actual neurons of the human brain.
Monotonicity of the activation function mirrors a type of biological neurons that fire
electrical pulses with a higher frequency if they receive stronger incoming signals. Bounded
functions are employed since biochemical limitations impose an upper bound on the firing frequency
of the biological neurons \cite{activation_glorot}. In modern
feed-forward neural networks, however, the \emph{Rectified Linear Unit (ReLU)} $\psi(z)=\max(z,0)$
is the most widely-used activation function \cite{relu}. The simple functional form
of the ReLU is advantageous not only from a design principle of minimalism, but also for
differentiating the whole neural network with respect to its parameters in the search for an optimal
configuration: Computing the derivative of the ReLU is trivial and the vanishing gradient problem in
the near constant regions of sigmoidal activation functions is avoided in the active region of the ReLU, which prevents the
gradient-based optimization methods in the subsequent paragraph from stalling \cite{activation_functions}. A more thorough discussion of activation
functions can be found in \cite{goodfellow}.

There is a downside of ANNs in comparison to a simple
linear regression: Firstly, it is no longer possible to find the optimal model parameters by
analytic means. Therefore, computationally efficient numerical optimization methods, such as
stochastic gradient descent, must be employed in the search for sufficiently
performant approximate solutions \cite{goodfellow}. Although these algorithms offer a practical
alternative to analytic solution strategies for
ANNs of reasonable size, the need for numerical optimization adds another layer of
complexity to the task of constructing regression models. Secondly, an ANN with enough parameters can always arrive at very low errors by
memorizing overly specific characteristics of the given data~\cite{overfitting}. As a remedy, the
data is split into a \emph{training set} and a \emph{test set}. The parameters of the ANN are
optimized with respect to the training set, which is called \emph{training}, and the performance of
the trained model is evaluated only with respect to the data in the test set. If the ANN performs
well on the test set, we can assume that it has learned a general rule that tells us how the input
variable is related to the target variable. Although the multitude of parameters makes it hard
for humans to understand this general rule, we know that it exists. In that sense, a low error on
the test set for an ANN gives us the same type of evidence for the existence of a relationship
between two entities, as a visual pattern or the Pearson correlation do in the $1$-dimensional case.

\section{Knot Invariants}
The techniques from the previous section were employed to probe possible relationships between several knot invariants. The investigated knot invariants can be divided into two groups: a group of combinatorial invariants is derived from the knot diagram, whereas a group of geometric invariants is derived from the hyperbolic structure on the knot complement. For the remainder of the article, let $K$ denote a hyperbolic knot in $S^3$.

\subsection{Combinatorial Knot Invariants}
The considered combinatorial invariants in the first group can be understood as compressed or generalized versions of the well-known Jones polynomial (see \cite{lickorish_jones} for the definition):
\begin{itemize}
	\item The Jones polynomial of $K$ is denoted by $J_K$.
	\item The \emph{knot determinant} of $K$ is given by $\det(K) = J_K(-1)$.
	\item The \emph{Mahler measure} of $J_K$ is given by $m(J_K) =  \exp\left(\frac{1}{2 \pi} \int_{0}^{2\pi} \ln\left|J_K\left(e^{i\theta}\right)\right| d\theta\right)$.
	\item As a slight variation of the knot determinant, we will also evaluate the Jones polynomial at the root of unity $\zeta = e^{2\pi i 3/5}$.
	\item The \emph{Khovanov polynomial} of $K$ is denoted by $\mathcal{KH}_K$. It is a Laurent polynomial in two variables whose graded Euler characteristic equals the Jones polynomial, in the sense that setting one of the variables to $-1$ recovers $J_K$. Mikhail Khovanov introduced it as a more refined version of the Jones polynomial \cite{khovanov_original}. For alternating knots, the diagonal coefficients equal the coefficients of the Jones polynomial and the off-diagonal coefficients vanish.
\end{itemize}

\subsection{Hyperbolic Knot Invariants}
The second group of knot invariants is derived by exploiting the geometric structure of the knot complement.
\begin{itemize}
	\item A \emph{hyperbolic knot} is a knot whose complement $S^3\setminus K$ admits a hyperbolic structure. Mostow rigidity \cite{mostow} shows that all such structures have the same volume and it is therefore a topological invariant, which is called the \emph{hyperbolic volume} $\vol{K}$.
	\item Witten \cite{witten} describes how the Jones polynomial can be defined using only intrinsic $3$\nobreakdash-dimensional properties of a knot. In this context, the Chern-Simons form from Chern-Simons theory can be used to assign a complex number to each hyperbolic knot whose modulus coincides up to a normalization factor with the hyperbolic volume. The argument of the complex number is called \emph{Chern-Simons invariant} $\mathcal{CS}(K) \in \mathbb{R} / \mathbb{Z}$ and it follows again from Mostow rigidity that it is a topological invariant of a hyperbolic knot.
	\item Another group of hyperbolic knot invariants is presented by Adams in \cite{adams}. Hyperbolic manifolds are seen as quotients $\mathbb{H}^3 / \Gamma$ of the hyperbolic $3$-space by a discrete group of fixed point free isometries $\Gamma$. If one considers tubular neighborhoods $T \times (0,t]$ of the deleted knot $T \times \{0\}$ in the knot complement, small neighborhoods lift to a disjoint set of horoballs in the covering space $\mathbb{H}^3$. Increasing the tubular neighborhood of the knot until the first two horoballs touch in the covering space, we obtain a maximal tubular neighborhood $T \times (0,t_{\max}]$ of the missing knot in the knot complement. This maximal tubular neighborhood is called the \emph{maximal cusp} and several hyperbolic knot invariants can be derived from it:
	\begin{itemize}
		\item The \emph{maximal cusp volume} $\mathcal{C}(K)$ is the hyperbolic volume of the maximal cusp.
		\item The \emph{longitude length} $\vert \lambda(K) \vert$ is the minimal length of a longitude in the boundary of the maximal cusp.
		\item The \emph{meridian length} $\vert \mu(K) \vert$ is the minimal length of a meridian in the boundary of the maximal cusp.
		\item The boundary torus of the maximal cusp has $\mathbb{R}^2$ as its universal cover and longitudes and meridians are non-trivial elements of the fundamental group $\pi_1(T \times \{t_{max}\})$. Therefore, we can go around the minimal length meridian and longitude once and observe the difference between the starting point and end point of a corresponding path in the cover $\mathbb{R}^2$. The respective differences are called \emph{longitude translation} $(\lambda_x, \lambda_y)$ and \emph{meridian translation} $(\mu_x, \mu_y)$. In order to compare the data for several knots, the orientations can be chosen such that we have $\lambda_y = 0$, which means that $\lambda_x$ is fully determined by $\vert \lambda(K) \vert$.
	\end{itemize}
\end{itemize}

\subsection{Known Connections}
To the best of our knowledge, theoretical research in knot theory has struggled to establish
connections between the combinatorial and geometric properties of knots. An exception is given by
the above-mentioned work of Witten that connects the Jones polynomial to the hyperbolic volume and
the Chern-Simons invariant via Chern-Simons theory \cite{witten}. Another example of speculative
nature is the long standing \emph{volume conjecture}, which predicts that the N-colored Jones
polynomials determine the volume of a hyperbolic knot \cite{kashaev, murakami}.
This conjecture is one of the most well-known and mysterious conjectures in low-dimensional topology
and it has been extensively studied. Its connection to Chern-Simons theory has been established in~\cite{murakami2002kashaevs}. We do not know of any conjecture relating the invariants derived from
the maximal cusp to the presented combinatorial ones.

Data-driven approaches have recently started to renew the interest in connections between
combinatorial and hyperbolic knot invariants:
Jejjala, Kar and Parrikar \cite{dl_jejjala} reported that a standard version of an artificial neural
network could predict the hyperbolic volume of a knot from the coefficients of
its Jones polynomial with a mean absolute percentage error of under 3\%. Their results for knots with up to 15~crossings
compared favorably to prediction models based on volume-ish
bounds and the Khovanov homology rank, respectively, and were robust to shrinking the size of the
training set. A spectral analysis of the weight matrices revealed a regularity in the largest
eigenvalues over several training runs, suggesting that these matrices encode a stable relationship
between the Jones polynomial and the hyperbolic volume. The authors proposed incorporating invariants from Khovanov homology for a possible
improvement of the predictive power, which is why we have included this aspect in our research.
A subsequent study conducted by Craven, Jejjala and Kar~\cite{craven} focussed on disentangling the approximate prediction formula encoded in the ANN from
\cite{dl_jejjala}. This was a necessary step for a further understanding of the precise relationship,
given that the original ANN was using over 10000 trainable parameters. In order to extract a more
compact and interpretable version of the formula that showed similar prediction performance,
layer-wise relevance propagation was employed to identify the features that were mainly responsible
for the successful prediction. They were able to distill a performant approximate formula of the
form
\begin{equation}
	\vol{K} \approx a \log(|J_K(e^{3\pi i / 4})| + b) - c,
\end{equation}
for real numbers $a, b, c$, which predicts the hyperbolic volume from a single evaluation
of the Jones polynomial at a root of unity. The simplicity of the formula is remarkable when
compared to the original formulation of the volume conjecture, which involves a limit over all the
N-colored Jones polynomials.
The success of these data-driven methods leads us to the the next section where we apply similar
techniques to analyze further connections between combinatorial and hyperbolic knot invariants.

\section{Experiments}
In this section, the conducted experiments and their results are documented and discussed.

\subsection{Data Preparation}
All of the experiments could not have been performed without the helpful software packages provided by the computational knot theory community. The following sources were used:

\begin{itemize}
	\item The Jones polynomials of all hyperbolic knots with up to 14 crossings (almost 60000 knots) were computed by a combination of SnapPy \cite{SnapPy} and SageMath \cite{sage}. The $1$-dimensional derived invariants were computed from the Jones polynomials with NumPy \cite{numpy}.
	\item The Khovanov polynomials for all hyperbolic knots with up to 12 crossings (almost 3000 knots) were computed with the Software package KhoHo \cite{khoho}.
	\item The hyperbolic volumes of all hyperbolic knots with up to 14 crossings were computed by SnapPy.
	\item Further hyperbolic knot invariants, namely the longitude length, the meridian length and translation, the maximal cusp volume and the Chern-Simons invariant, were taken from the online database KnotInfo \cite{knotinfo}. The database includes several knot invariants for all prime knots with up to 12 crossings.
\end{itemize}

The coefficients of the Jones polynomials and the Khovanov polynomials were vectorized and padded with zeros to obtain a consistent data format. The three single-variable invariants $\det(K), m(J_K)$ and $J_K(\zeta)$ were transformed to a logarithmic scale and divided by $\ln(\deg(J_K)))$, in order to obtain the rescaled invariants $\underline{\det}(K),\underline{m}(J_K)$ and $\underline{J}_K(\zeta)$.\footnote{The \emph{degree} $\deg(J_K)$ of a Jones polynomial $J_K(t) = c_0t^{n_0} + ... + c_kt^{n_k}$ with $c_0 \neq 0$ and $c_k \neq 0$ is given by $n_k - n_0.$} All hyperbolic invariants were left unchanged.

\subsection{Linear Regression}
The first part of our experiments focussed on using linear regression to detect possible connections
between the $1$-dimensional knot invariants.
\subsubsection{Model Design}
Separate models were constructed for alternating knots $\mathcal{K}_a$, non-alternating knots $\mathcal{K}_n$ and all knots $\mathcal{K}$. For each pair $(\Phi,\Psi)$ with
\begin{equation}
	\Phi \in \{\underline{\det}(K),\underline{m}(J_K),\underline{J}_K(\zeta)\}
\end{equation}
and
\begin{equation}
	\Psi \in \{ \vol{K},\vert\lambda(K)\vert,\vert\mu(K)\vert,\mu_x(K),\mu_y(K),\mathcal{C}(K),\mathcal{CS}(K) \},
\end{equation}
a least-squares linear regression was fit using NumPy and the Pearson correlation was computed.

\subsubsection{Experimental Results}

\begin{table}[htb!]
	\begin{tabular}{l|lll}
		\multicolumn{1}{c|}{}  & \multicolumn{3}{c}{\boldmath $r$\unboldmath}\\
		\hline
		& $\mathcal{K}$ & $\mathcal{K}_{a}$ & $\mathcal{K}_{n}$\\
		\hline
		$\underline{\det}(K)$ & 0.74 & 0.98 & 0.66\\
		$\underline{m}(J_K)$ & 0.78 & 0.71 & 0.77\\
		$\underline{J}_K(\zeta)$ & 0.94 & 0.96 & 0.92\\
	\end{tabular}
	\caption{Pearson correlation coefficients between three combinatorial knot invariants and the hyperbolic volume for knots with up to 14 crossings.\label{corr_14}}
\end{table}

\begin{figure} [htb!]
	\begin{center}
		\includegraphics[width=12cm]{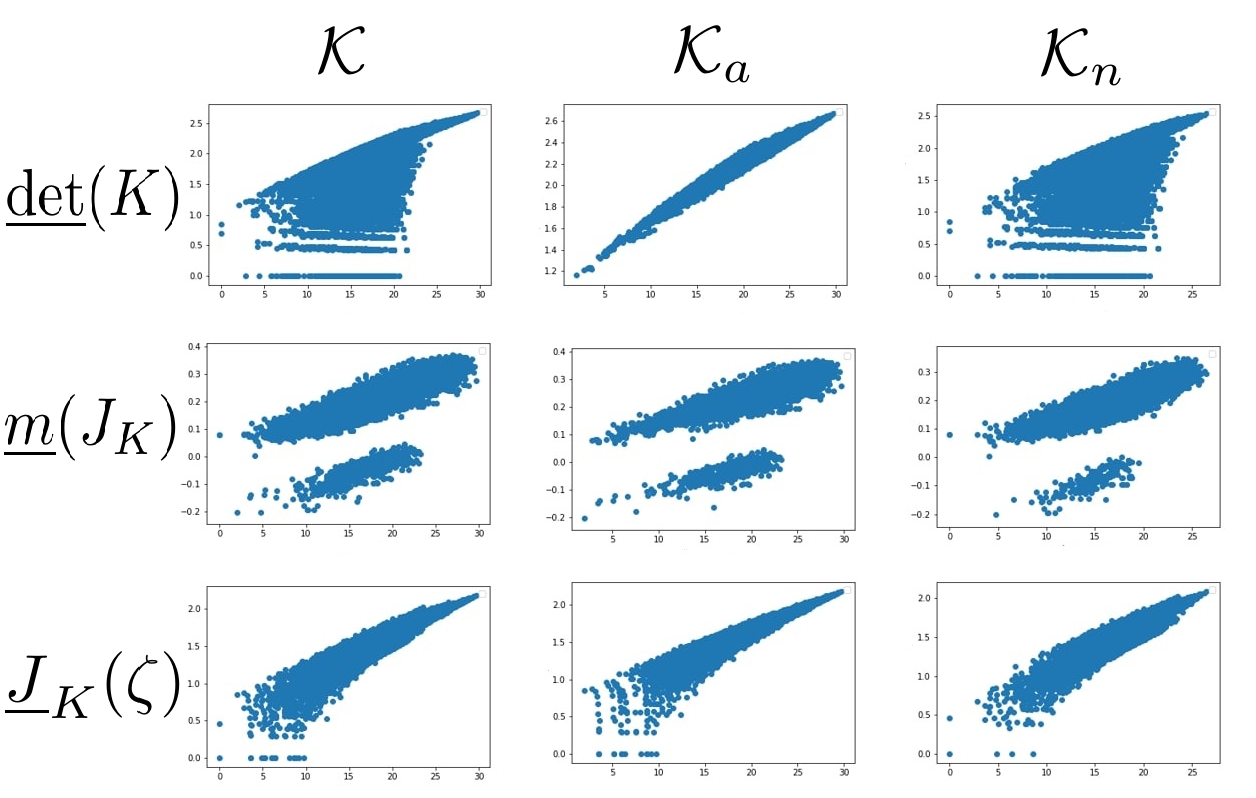}
		\caption[crossval]{Plots of the hyperbolic volumes corresponding to the linear regression evaluations from Table \ref{corr_14}.}
		\label{plot14_combined}
	\end{center}
\end{figure}

For hyperbolic knots with up to 14 crossings, the approximately linear relationships between $\vol{K}$ on the one hand and $\underline{\det}(K),\underline{m}(J_K)$ and $\underline{J}_K(\zeta)$ on the other are shown qualitatively in Fig. \ref{plot14_combined} and quantitatively in Table~\ref{corr_14}.

\begin{itemize}
	\item The first row shows how the relationship between the knot determinant and the hyperbolic volume varies drastically depending on the knots being alternating or not ($r=0.98$ vs. $r=0.66$), as was already observed by Dunfield \cite{dunfield_determinant} and Friedl and Jackson \cite{friedl_hyperbolic}. The correlation on all knots is $r=0.74$.
	\item The second row shows the more stable behavior for the Mahler measure (${r=0.71}$~vs.~$r=0.77$) and the characteristic formation of two clusters. Given that a linear regression fits only one line through the data, the overall correlation of $r=0.77$ barely surpasses the knot determinant. Fitting two separate lines through the data from all knots yields intracluster correlations of $r=0.87$ and $0.82$, respectively. The separation into two clusters did not occur for the Mahler measure of the Alexander polynomial in \cite{friedl_hyperbolic}.
	\item The third row shows that the evaluation of
the Jones polynomial at the complex root of unity performs well in both cases ($r= 0.96$ vs.  $r=0.92$) without splitting the knots into two groups, leading to the best overall performance of $r=0.94$. The correlation is still weaker than for the twisted Alexander polynomial that was proposed by Friedl and Jackson ($r=0.98$ for knots with up to 15 crossings), but the simplicity of evaluating the Jones polynomial is remarkable in comparison from a theoretical and computational point of view.
\end{itemize}

Table \ref{corr_12} shows the correlation coefficients between combinatorial and hyperbolic
invariants for hyperbolic knots with up to 12 crossings. It is not surprising that the first column again
shows high values, given that the used knot data is a subset of the knot data that was used to create Table \ref{corr_14}. Furthermore, it is important to keep in mind that being linearly related is an equivalence relation, meaning that high correlations between the combinatorial invariants and hyperbolic invariants other than $\vol{K}$ already imply high correlations between the hyperbolic invariant and $\vol{K}$. This is most prominent for the maximal cusp volume $\mathcal{C}(K)$, whose high correlation with $\vol{K}$ could be confirmed by plotting both invariants against each other. However, observing these implied correlations with $\vol{K}$ can be interesting in its own right.

\begin{table}[htb!]
	\begin{tabular}{l|lllllll}
			\hline
			\hline
			\multicolumn{1}{c|}{}  & \multicolumn{7}{c}{\boldmath $r$ for $\mathcal{K}$ \unboldmath}\\
			\hline
			& $\vol{K}$ & $\vert\lambda(K)\vert$ & $\vert\mu(K)\vert$ & $\mu_x(K)$ & $\mu_y(K)$ & $\mathcal{C}(K)$ & $\mathcal{CS}(K)$\\
			\hline
			$\underline{\det}(K)$ & 0.80 & 0.26 & -0.41 & -0.38 & 0.27 & 0.76 & -0.03\\
			$\underline{m}(J_K)$ & 0.70 & 0.27 & -0.17 & -0.22 & 0.22 & 0.66 & -0.02\\
			$\underline{J}_K(\zeta)$ & 0.94 & 0.25 & -0.34 & -0.44 & 0.39 & 0.90 & -0.03\\
			\hline
			\hline
			\multicolumn{1}{c|}{}  & \multicolumn{7}{c}{\boldmath $r$ for $\mathcal{K}_a$ \unboldmath}\\
			\hline
			& $\vol{K}$ & $\vert\lambda(K)\vert$ & $\vert\mu(K)\vert$ & $\mu_x(K)$ & $\mu_y(K)$ & $\mathcal{C}(K)$ & $\mathcal{CS}(K)$\\
			\hline
			$\underline{\det}(K)$ & 0.99 & 0.28 & 0.02 & -0.37 & 0.47 & 0.91 & -0.05\\
			$\underline{m}(J_K)$ & 0.65 & 0.26 & 0.05 & -0.15 & 0.25 & 0.61 & -0.02\\
			$\underline{J}_K(\zeta)$ & 0.95 & 0.23 & -0.03 & -0.40 & 0.49 & 0.89 & -0.05\\
			\hline
			\hline
			\multicolumn{1}{c|}{}  & \multicolumn{7}{c}{\boldmath $r$ for $\mathcal{K}_n$ \unboldmath}\\
			\hline
			& $\vol{K}$ & $\vert\lambda(K)\vert$ & $\vert\mu(K)\vert$ & $\mu_x(K)$ & $\mu_y(K)$ & $\mathcal{C}(K)$ & $\mathcal{CS}(K)$\\
			\hline
			$\underline{\det}(K)$ & 0.73 & 0.15 & -0.47 & -0.35 & 0.23 & 0.71 & 0.00\\
			$\underline{m}(J_K)$ & 0.73 & 0.20 & -0.31 & -0.25 & 0.21 & 0.70 & 0.00\\
			$\underline{J}_K(\zeta)$ & 0.93 & 0.16 & -0.46 & -0.42 & 0.34 & 0.88 & 0.00\\
	\end{tabular}
	\caption{Pearson correlation coefficients between combinatorial and hyperbolic invariants for knots with up to 12 crossings. \label{corr_12}}
\end{table}

\begin{itemize}
	\item The Chern-Simons invariant is the only invariant that shows only very low values. In light of the already established theoretical connection between the Jones polynomial and the Chern-Simons invariant via Chern-Simons theory, the negative results are somewhat surprising. The least-squares optimization might have been inappropriate in this case, given that the invariant takes values in $\mathbb{R} / \frac{1}{2}\mathbb{Z}$ for the data from KnotInfo and we used the Euclidean distance between the representatives without taking into account the cyclic structure.
	\item All of the other invariants show moderate connections to $\vol{K}$, with absolute values of most coefficients ranging from $0.2$ to $0.5$.
	\item Notable exceptions are the meridian length $\vert\mu(K)\vert$, where the correlation coefficient nearly vanishes for alternating knots, and the longitude length $\vert\lambda(K)\vert$, where the correlation weakens for non-alternating knots.
\end{itemize}

\subsection{Deep Learning}
For the multidimensional input invariants, we replaced the linear regression models by fully
connected ANNs and proceeded similarly to the previously presented experiments.
\subsubsection{Model Design}
As in the linear case above, separate models were constructed for alternating knots~$\mathcal{K}_a$, non-alternating knots $\mathcal{K}_n$ and all knots $\mathcal{K}$.
\begin{itemize}
	\item For each pair of invariants $(\Phi,\Psi)$ with input
	\begin{equation}
		\Phi \in \{J_K, \mathcal{KH}_K\}
	\end{equation}
	and target
	\begin{equation}
		\Psi \in \{ \vol{K},\vert\lambda(K)\vert,\vert\mu(K)\vert,\mu_x(K),\mu_y(K),\mathcal{C}(K),\mathcal{CS}(K) \},
	\end{equation}
	 a fully connected ANN with two hidden layers, consisting of 100 hidden neurons each and using
	 ReLU activation functions, was set up in the deep learning framework TensorFlow \cite{tensorflow}. The
	 model was trained on randomly selected 80\% of the data and evaluated on the remaining 20\%,
	 with respect to the prediction errors $\Delta_{\text{MAPE}}$ and $\Delta_{\text{MSE}}$.
	\item For each pair of invariants $(\Phi,\Psi)$ with input
	\begin{equation}
		\Phi \in \{\underline{\det}(K),\underline{m}(J_K),\underline{J}_K(\zeta)\}
	\end{equation}
	 and target $\Psi$ as above, a least-squares linear regression was fit on randomly selected 80\% of the data. Subsequently, the prediction errors $\Delta_{\text{MAPE}}$ and $\Delta_{\text{MSE}}$ were computed on the remaining 20\% of the data, in order to compare the predictive power of the linear models and the ANNs.
\end{itemize}

\subsubsection{Experimental Results}
The prediction errors of the ANNs and linear regressions are reported in Table~\ref{mape2_12} for
mean absolute percentage errors and Table \ref{mse2_12} for mean squared errors, respectively. The
last row of the tables shows the errors of a model that always predicts the mean value of the
hyperbolic target invariant. This serves as a base line that allows a comparison of the predictive
power of the optimized statistical models and a simple educated guess. Prediction errors are grouped
by horizontal lines according to the model type used (ANN, linear regression or base line model).
The mean squared errors are reported as fractions of the mean squared errors of the base line model,
in order to facilitate comparison.

\begin{table}[htb!]
	\begin{tabular}{l||lllllll}
			\hline
			\hline
			\hline
			\multicolumn{1}{c||}{}  & \multicolumn{7}{c}{\boldmath $\Delta_{\text{MAPE}}$ for $\mathcal{K}$ \unboldmath}\\
			\hline
			\hline
			& $\vol{K}$ & $\vert\lambda(K)\vert$ & $\vert\mu(K)\vert$ & $\mu_x(K)$ & $\mu_y(K)$ & $\mathcal{C}(K)$ & $\mathcal{CS}(K)$\\
			\hline
			\hline
			$\mathcal{KH}_K$ & \bf{12.1} & \bf{8.9} & 6.6 & - & \bf{12.1} & 14.7 & -\\
			$J_K$ & \bf{3.9} & 15.0 & 6.3 & - & \bf{11.7} & \bf{8.8} & -\\
			\hline
			$\underline{\det}(K)$ & 13.4 & 23.4 & 6.7 & - & 22.8 & 14.9 & -\\
			$\underline{m}(J_K)$ & 15.1 & 23.4 & 7.1 & - & 23.1 & 16.7 & -\\
			$\underline{J}_K(\zeta)$ & \bf{7.2} & 22.9 & 6.9 & - & 21.4 & \bf{9.5} & -\\
			\hline
			base line & 24.2 & 25.7 & 7.1 & - & 24.6 & 24.2 & -\\
			\hline
			\hline
			\hline
			\multicolumn{1}{c||}{}  & \multicolumn{7}{c}{\boldmath $\Delta_{\text{MAPE}}$ for $\mathcal{K}_a$ \unboldmath}\\
			\hline
			\hline
			& $\vol{K}$ & $\vert\lambda(K)\vert$ & $\vert\mu(K)\vert$ & $\mu_x(K)$ & $\mu_y(K)$ & $\mathcal{C}(K)$ & $\mathcal{CS}(K)$\\
			\hline
			\hline
			$\mathcal{KH}_K$ & \bf{5.8} & \bf{5.5} & 6.0 & - & \bf{7.3} & 10.0 & -\\
			$J_K$ & \bf{2.7} & \bf{5.9} & 5.9 & - & \bf{8.3} & \bf{7.1} & -\\
			\hline
			$\underline{\det}(K)$ & \bf{3.0} & 17.5 & 6.5 & - & 19.4 & \bf{7.9} & -\\
			$\underline{m}(J_K)$ & 14.3 & 18.0 & 6.5 & - & 21.5 & 14.2 & -\\
			$\underline{J}_K(\zeta)$ & \bf{5.6} & 18.0 & 6.5 & - & 18.7 & \bf{8.3} & -\\
			\hline
			base line & 20.0 & 19.1 & 6.5 & - & 22.8 & 19.1 & -\\
			\hline
			\hline
			\hline
			\multicolumn{1}{c||}{}  & \multicolumn{7}{c}{\boldmath $\Delta_{\text{MAPE}}$ for $\mathcal{K}_n$ \unboldmath}\\
			\hline
			\hline
			& $\vol{K}$ & $\vert\lambda(K)\vert$ & $\vert\mu(K)\vert$ & $\mu_x(K)$ & $\mu_y(K)$ & $\mathcal{C}(K)$ & $\mathcal{CS}(K)$\\
			\hline
			\hline
			$\mathcal{KH}_K$ & 12.5 & \bf{13.5} & 6.9 & - & \bf{13.0} & 15.9 & -\\
			$J_K$ & \bf{3.9} & \bf{15.0} & 7.0 & - & 14.0 & \bf{9.2} & -\\
			\hline
			$\underline{\det}(K)$ & 13.8 & 28.9 & 5.9 & - & 25.9 & 16.1 & -\\
			$\underline{m}(J_K)$ & 13.2 & 28.2 & 6.4 & - & 26.4 & 14.6 & -\\
			$\underline{J}_K(\zeta)$ & \bf{6.9} & 28.8 & 6.0 & - & 23.8 & \bf{9.3} & -\\
			\hline
			base line & 20.7 & 30.5 & 6.4 & - & 27.6 & 22.7 & -\\
	\end{tabular}
	\caption{Mean absolute percentage errors for predictions of hyperbolic knot invariants of knots with up to 12 crossings. Within each subtable, errors are grouped by horizontal lines according to the model type used (ANN, linear regression or base line model) and printed in bold if they are less than half of the error of the base line model. Two columns could not be calculated because the data contained zeros.\label{mape2_12}}
\end{table}

\begin{table}[htb!]
	\begin{tabular}{l||lllllll}
			\hline
			\hline
			\hline
			\multicolumn{1}{c||}{}  & \multicolumn{7}{c}{\boldmath Relative $\Delta_{\text{MSE}}$ for $\mathcal{K}$ \unboldmath}\\
			\hline
			\hline
			& $\vol{K}$ & $\vert\lambda(K)\vert$ & $\vert\mu(K)\vert$ & $\mu_x(K)$ & $\mu_y(K)$ & $\mathcal{C}(K)$ & $\mathcal{CS}(K)$\\
			\hline
			\hline
			$\mathcal{KH}_K$ & \bf{0.31} & \bf{0.24} & 1.00 & \bf{0.22} & \bf{0.34} & \bf{0.46} & 1.19\\
			$J_K$ & \bf{0.03} & \bf{0.39} & 0.88 & \bf{0.36} & \bf{0.46} & \bf{0.19} & 0.92\\
			\hline
			$\underline{\det}(K)$ & \bf{0.35} & 0.89 & 0.85 & 0.87 & 0.94 & \bf{0.43} & 1.00\\
			$\underline{m}(J_K)$ & \bf{0.47} & 0.92 & 0.96 & 0.94 & 0.93 & 0.55 & 1.00\\
			$\underline{J}_K(\zeta)$ & \bf{0.10} & 0.89 & 0.89 & 0.83 & 0.88 & \bf{0.20} & 1.00\\
			\hline
			base line & 1.00 & 1.00 & 1.00 & 1.00 & 1.00 & 1.00 & 1.00\\
			\hline
			\hline
			\hline
			\multicolumn{1}{c||}{}  & \multicolumn{7}{c}{\boldmath Relative $\Delta_{\text{MSE}}$ for $\mathcal{K}_a$ \unboldmath}\\
			\hline
			\hline
			& $\vol{K}$ & $\vert\lambda(K)\vert$ & $\vert\mu(K)\vert$ & $\mu_x(K)$ & $\mu_y(K)$ & $\mathcal{C}(K)$ & $\mathcal{CS}(K)$\\
			\hline
			\hline
			$\mathcal{KH}_K$ & \bf{0.13} & \bf{0.10} & 0.96 & \bf{0.08} & \bf{0.17} & \bf{0.32} & 1.01\\
			$J_K$ & \bf{0.02} & \bf{0.14} & 1.10 & \bf{0.10} & \bf{0.20} & \bf{0.18} & 0.55\\
			\hline
			$\underline{\det}(K)$ & \bf{0.03} & 0.96 & 1.00 & 0.88 & 0.80 & \bf{0.20} & 0.99\\
			$\underline{m}(J_K)$ & 0.55 & 0.96 & 0.99 & 0.98 & 0.92 & 0.59 & 1.00\\
			$\underline{J}_K(\zeta)$ & \bf{0.09} & 0.99 & 0.99 & 0.85 & 0.76 & \bf{0.21} & 0.99\\
			\hline
			base line & 1.00 & 1.00 & 1.00 & 1.00 & 1.00 & 1.00 & 1.00\\
			\hline
			\hline
			\hline
			\multicolumn{1}{c||}{}  & \multicolumn{7}{c}{\boldmath Relative $\Delta_{\text{MSE}}$ for $\mathcal{K}_n$ \unboldmath}\\
			\hline
			\hline
			& $\vol{K}$ & $\vert\lambda(K)\vert$ & $\vert\mu(K)\vert$ & $\mu_x(K)$ & $\mu_y(K)$ & $\mathcal{C}(K)$ & $\mathcal{CS}(K)$\\
			\hline
			\hline
			$\mathcal{KH}_K$ & \bf{0.45} & \bf{0.35} & 1.12 & \bf{0.28} & \bf{0.37} & 0.59 & 1.13\\
			$J_K$ & \bf{0.05} & \bf{0.32} & 1.09 & \bf{0.46} & 0.56 & \bf{0.22} & 1.16\\
			\hline
			$\underline{\det}(K)$ & 0.52 & 0.94 & 0.71 & 0.89 & 1.05 & 0.54 & 1.01\\
			$\underline{m}(J_K)$ & 0.55 & 0.91 & 0.89 & 0.94 & 1.01 & 0.61 & 1.00\\
			$\underline{J}_K(\zeta)$ & \bf{0.13} & 0.95 & 0.80 & 0.81 & 0.96 & \bf{0.19} & 1.00\\
			\hline
			base line & 1.00 & 1.00 & 1.00 & 1.00 & 1.00 & 1.00 & 1.00\\
	\end{tabular}
	\caption{Relative mean squared errors for predictions of hyperbolic knot invariants of knots with up to 12 crossings. Within each subtable, errors are grouped by horizontal lines according to the model type used (ANN, linear regression or base line model) and printed in bold if they are less than half of the error of the base line model.\label{mse2_12}}
\end{table}

\begin{itemize}
	\item As was already shown by Jejjala, Kar and Parrikar in \cite{dl_jejjala}, the hyperbolic volume can be predicted with high precision by an ANN that takes the coefficients of the Jones polynomial as an input. Unfortunately, using the coefficients of the Khovanov polynomials did not improve the predictive power of the model. On the contrary, the observed performance decreased considerably, even though all the information from the Jones polynomials is contained in the Khovanov polynomials. This may indicate that the additional information in the Khovanov polynomials does not add any valuable input to the prediction. Since the coefficients of the Khovanov polynomials are two-dimensional matrices, the dimension of the input increases drastically when compared to the prediction based on the Jones polynomials. This makes the optimization problem in the training process harder, as more parameters in the input layer need to be optimized, and it can lead to a worse resulting model. Further evidence for this hypothesis is given by the fact that the difference between the prediction errors of the two models was smaller for alternating knots, where most of the Khovanov coefficients are zero. The general statement that high-dimensional problems require more data for accurate predictions is known as the \emph{curse of dimensionality} \cite{curse}.
	\item For the prediction of the longitude length, the situation changed and the Khovanov ANN was the best model by a large margin, even though the Jones ANN was also clearly better than the base line.
	\item None of the models could predict the meridian length significantly better than the base line.
	\item For the meridian translation in both directions, the Khovanov ANNs were again the best predictors, followed by the Jones ANNs, which also clearly beat the linear models and the base line.
	\item The maximal cusp volume shows similar results to the hyperbolic volume, which is not surprising in light of the above discussion about their strong linear association. It is worth noting that the Jones ANN performed considerably worse than in the case of predicting the hyperbolic volume, which could be due to optimization problems in the network.
	\item The predictions for the Chern-Simons invariant were the least accurate ones, with only one small success: In the case of alternating knots, the Jones ANN beat the base line by a large margin. As discussed above, the poor results for the other input variables could be due to the neglected cyclic nature of the Chern-Simons invariant.
\end{itemize}

\subsection{Varying Regression Inputs and Network Size}
Apart from evaluating the Jones polynomial at $\zeta = e^{2\pi i 3/5}$, other roots of unity of the form $e^{2\pi i k/n}$ were tested. It turned out that choosing $k/n$ close, but not equal to $1/2$ provided the best results, and $k=3$, $n=5$ was chosen for simplicity reasons.
Inspired by the success of $\underline{J}_K(\zeta)$, ANNs were tested that input the values of $J_K$ at regularly spaced roots of unity instead of the coefficients of $J_K$. The results did not differ significantly.
Most experiments were repeated with smaller ANNs in order to probe the necessary representative
power of the models. The hyperbolic volume could be predicted with under 5\% error from an ANN with
only one hidden layer and 5 hidden neurons, reducing the number of parameters in the weight matrices
from over 10000 to only 80. In principle, it is possible to read off a formula relating the Jones
polynomial to the hyperbolic volume with high accuracy from these smaller networks, but the number
of parameters is still far too high to allow a conceptual interpretation of the resulting
expression. As mentioned above, Craven, Jejjala and Kar \cite{craven} have recently solved this very
problem by using layer-wise relevance propagation and distilled the compact formula
\begin{equation} \label{formula_craven}
	\vol{K} \approx 6.20 \log(|J_K(e^{3\pi i / 4})| + 6.77) - 0.94,
\end{equation}
which needs only one evaluation of the Jones polynomial at a root of unity to predict the hyperbolic
volume with an error of 2.86\% on all hyperbolic knots up to and including 16 crossings. Note the
similarity of this formula to the linear regression using $\underline{J}_K(\zeta)$ that was
discussed above. Our way of reasoning was to replace the root of unity $-1$ in the definition of
the knot determinant by another root of unity, whereas (\ref{formula_craven}) was extracted
in a purely data-driven way from the complex inner
workings of an ANN with several thousands of parameters. Given that the numerical
constants in (\ref{formula_craven}) are derived from a least-squares optimization on a finite dataset, it is natural that it
is not an exact formula. Apart from the trainable parameters, the phase in the exponential is a
hyperparameter that can be varied to some degree without losing much predictive power, which is
investigated in more detail in \cite{craven}. These observations raise the question how the above
approximate formula can be turned into a theorem with meaningful and interpretable constants. Do we
need to include additional variables, as was the case for Davies et al. in \cite{davies}? Are there
multiple phases that lead to different exact formulas or is there a single viable choice? Similar
questions can be asked about the empirical findings in our own work, which leads us to the last
section.

\section{Open Questions}
In summary, we have seen that statistical methods, such as linear regression and deep learning with artificial neural networks, can be valuable tools for the discovery of relationships between knot invariants. Now that we know of some empirical connections between combinatorial and hyperbolic knot invariants, further mathematical research could address the following questions:
\begin{itemize}
	\item How can the strong correlation between $\underline{J}_K(\zeta)$ and $\vol{K}$ be explained?
	\item Why is $\underline{\det}(K)$ only effective for predicting $\vol{K}$ in the case of alternating knots?
	\item What is a theoretical reason for the formation of the two clusters when comparing $\underline{m}(J_K)$ to $\vol{K}$ in Fig. \ref{plot14_combined}?
	\item What is an exact formula that predicts $\vol{K}$ from $J_K$? Are additional invariants
	needed for a perfect prediction?
	\item What is the theoretical relationship between $\mathcal{KH}_K$ and $\mu(K)$? Are additional invariants needed for an exact prediction?
	\item Why is there a clear correlation between $\vol{K}$ and $|\mu(K)|$ only for non-alternating knots?
	\item How can $\mathcal{CS}(K)$ be related to $J_K$ in the alternating case? What other invariants are needed for an exact formula? Is there a relationship in the non-alternating case?
	\item What happens if we repeat the experiments with links instead of knots?
\end{itemize}
Clearly, most of these questions cannot be answered
by the presented statistical methods alone and need to be tackled by mathematical research of a more traditional nature.
ANNs are simply not made to be interpretable, as was already discussed by Hughes
\cite{hughes}: they are designed to achieve a maximum of predictive power based on many non-linear
transformations of the input variables. Perhaps the growing interest in explainable artificial
intelligence \cite{explainable_ai} will produce algorithms that can support mathematicians even
further in the quest for detecting and deciphering relationships between mathematical entities.

\section*{Acknowledgements}
This article presents work that was conducted as part of my master thesis \emph{Deep Learning of Hyperbolic Knot Invariants} \cite{masterarbeit} at Universität Regensburg. The thesis provides a more detailed exposition of all the concepts that are discussed in this article and it can be accessed online. I would like to thank Stefan Friedl, the advisor of the thesis, for many helpful conversations and suggestions. I am also grateful to Lukas Lewark for his input on topics in computational knot theory. Maike Stern and Heribert Wankerl deserve special credit for proofreading the paper.

\bibliographystyle{unsrt}
\bibliography{references_paper}
\end{document}